\documentclass[12pt,reqno]{amsart}
\usepackage{amsmath,amssymb,amsfonts,amscd,latexsym,amsthm,mathrsfs,verbatim}
\usepackage[unicode]{hyperref}
\usepackage{fouridx}

\renewcommand{\d}{{\mathrm d}}
\newcommand{\mm}{{\mathrm m}}

\theoremstyle{remark}
\newtheorem*{question}{\bf Question}
\begin{document}

\title{Ap\'ery limits for elliptic $L$-values}

\author{Christoph Koutschan}
\address{
        Johann Radon Institute for Computational and Applied Mathematics (RICAM),
        Austrian Academy of Sciences,
        Altenberger Stra\ss e 69,
        A-4040 Linz, Austria}
\email{christoph.koutschan@ricam.oeaw.ac.at}

\author{Wadim Zudilin}
\address{
Department of Mathematics, IMAPP, Radboud University, PO Box 9010, 6500\,GL Nijmegen, Netherlands}
\email{w.zudilin@math.ru.nl}

\date{16 November 2021}

\begin{abstract}
For an (irreducible) recurrence equation with coefficients from $\mathbb Z[n]$ and its two linearly independent rational solutions $u_n,v_n$, the limit of $u_n/v_n$ as $n\to\infty$, when exists, is called the Ap\'ery limit.
We give a construction that realises certain quotients of $L$-values of elliptic curves as Ap\'ery limits.
\end{abstract}



\thanks{Research of the first author is supported by the Austrian Science Fund (FWF) grant F5011-N15.
Research of the second author is supported by the Dutch Research Council (NWO) grant OCENW.KLEIN.006.}

\maketitle

Ap\'ery's famous proof \cite{Po79} of the irrationality of $\zeta(3)$ displayed a particular phenomenon (which could have been certainly dismissed if discussed in the arithmetic context of some \emph{boring} quantities).
One considers the recurrence equation
\begin{equation}
(n+1)^3v_{n+1}-(2n+1)(17n^2+17n+5)v_n+n^3v_{n-1}=0
\quad\text{for}\; n=1,2,\dotsc.
\label{ap}
\end{equation}
and its two \emph{rational} solutions $u_n$ and $v_n$, where $n\ge0$, originating from the initial data $u_0=0$, $u_1=6$ and $v_0=1$, $v_1=5$.
Then $v_n$ are in fact integral for any $n\ge0$ and the denominators of $u_n$ have a moderate growth with $n$\,---\,certainly not like $n!^3$ as suggested by the recursion\,---\,but $O(C^n)$ as $n\to\infty$, for some $C>1$.
Namely, $D_n^3u_n\in\mathbb Z$ for all $n\ge1$, where $D_n$ denotes the least common multiple of $1,2,\dots,n$;
the asymptotics $D_n^{1/n}\to e$ as $n\to\infty$ is a consequence of the prime number theorem.
An important additional property is that the quotient $u_n/v_n\to\zeta(3)$ as $n\to\infty$ (and also $u_n/v_n\ne\zeta(3)$ for \emph{all}~$n$); even sharper: $v_n\zeta(3)-u_n\to0$ as $n\to\infty$; and at the highest level of sharpness we have $D_n^3(v_n\zeta(3)-u_n)\to0$ as $n\to\infty$.
It is the latter sharpest form that leads to the conclusion $\zeta(3)\notin\mathbb Q$.
But already the arithmetic properties of $u_n,v_n$ coupled with the `irrational' limit relation $u_n/v_n\to\zeta(3)$ as $n\to\infty$ are phenomenal.

One way to prove all the above claims in one shot is to cast the sequence $I_n=v_n\zeta(3)-u_n$ as the Beukers triple integral~\cite{Be79}
\[
I_n=\frac12\int_0^1\!\int_0^1\!\int_0^1\frac{x^n(1-x)^ny^n(1-y)^nz^n(1-z)^n}{(1-(1-xy)z)^{n+1}}\,\d x\,\d y\,\d z
\quad\text{for}\; n=1,2,\dotsc.
\]
A routine use of creative telescoping machinery, based on the Almkvist--Zeilberger algorithm \cite{AZ90} (in fact, its multivariable version \cite{AZ06}), then shows that $I_n$ indeed satisfies \eqref{ap}, while the evaluations $I_0=\zeta(3)$ and $I_1=5\zeta(3)-6$ are straightforward.
The arithmetic and analytic properties follow from the analysis of the integrals $I_n$ performed in \cite{Be79};
more \emph{practically}, they can be predicted/checked numerically based on the recurrence equation~\eqref{ap}.

A common belief is that we have a better understanding of the phenomenon these days.
Namely, we possess some (highly non-systematic!) recipes and strategies (see, for example, \cite{ASZ08,CS20,D-BKZ21,Za09,Zu21a,Zu21b}) for getting other meaningful constants $c$ as \emph{Ap\'ery limits}\,---\,in other words, there are (irreducible) recurrence equations with coefficients from $\mathbb Z[n]$ such that for two \emph{rational} solutions $u_n,v_n$ we have $u_n/v_n\to c$ as $n\to\infty$ and the denominators of $u_n,v_n$ are growing at most exponentially in~$n$.
(We may also consider \emph{weak} Ap\'ery limits when the latter condition on the growth of denominators is dropped.)
Though one would definitely like to draw some conclusions about the irrationality of those constants $c$, this constraint for the arithmetic to be in the sharpest form would severly shorten the existing list of known Ap\'ery limits;
for example, it would throw out Catalan's constant from the list.
A very basic question is then as follows.

\begin{question}
What real numbers can be realised as Ap\'ery limits?
\end{question}

Without going at any depth into this direction, we present here a (`weak') construction of Ap\'ery limits which are related to the $L$-values of elliptic curves (or of weight~2 modular forms).
The construction emanates from identities, most of which remain conjectural, between the $L$-values and Mahler measures.

\medskip
Consider the family of double integrals
\begin{align*}
J_n(z)
&=\int_0^1\!\int_0^1\frac{x^{n-1/2}(1-x)^{n-1/2}y^{n-1/2}(1-y)^n}{(1-zxy)^{n+1/2}}\,\d x\,\d y
\\
&=\frac{\Gamma(n+\tfrac12)^3\Gamma(n+1)}{\Gamma(2n+1)\Gamma(2n+\tfrac32)}
\cdot{}_3F_2\bigg(\begin{matrix} n+\tfrac12, \, n+\tfrac12, \, n+\tfrac12 \\ 2n+1, \, 2n+\tfrac32 \end{matrix} \biggm|z\bigg).
\end{align*}
Thanks to the nice hypergeometric representation, a recurrence equation
satisfied by the double integral can be computed using Zeilberger's fast
summation algorithm~\cite{AZ06,Ze90}, which is based on the method of creative
telescoping.  It leads to the following third-order recurrence equation:
\begin{align*}
&
4z^4(2n+1)^2(n+1)^2 \big(16(27z-32)n^4 - 16(69z-86)n^3
\\ &\quad
+ 8(108z-143)n^2 - 4(55z-76)n + 3(7z-10)\big) J_{n+1}
\displaybreak[2]\\ &\;
+ z^2\big(256(3z+8)(27z-32)n^8 - 256(3z+8)(15z-22)n^7
\\ &\;\quad
- 64(651z^2+661z-1744)n^6 + 192(59z^2-186)n^5
\\ &\;\quad
+ 16(1503z^2+697z-3610)n^4 - 16(79z^2-290z+116)n^3
\\ &\;\quad
- 4(569z^2-381z-580)n^2 + 4(11z^2-44z+18)n + 3(4z+3)(7z-10)\big) J_n
\displaybreak[2]\\ &\;
+ 4n\big(64(3z^2-20z+16)(27z-32)n^7 - 384(3z^2-20z+16)(7z-9)n^6
\\ &\;\quad
- 16(411z^3-2698z^2+3988z-1696)n^5 + 64(183z^3-1372z^2+2339z-1134)n^4
\\ &\;\quad
+ 4(531z^3-1400z^2-424z+1240)n^3 - 8(571z^3-4001z^2+6532z-3060)n^2
\\ &\;\quad
+ (151z^3-4742z^2+11596z-6888)n + 12(14z^2-29z-30)(z-1)\big) J_{n-1}
\displaybreak[2]\\ &\;
+ 4n(n-1)(2n-3)^2(z-1)\big(16(27z-32)n^4 + 48(13z-14)n^3
\\ &\;\quad
+ 8(18z-11)n^2 - 4(19z-24)n - (7z+6)\big) J_{n-2}
=0.
\end{align*}
Furthermore, if we take
\begin{align*}
\lambda(z)
&=J_0(z)=2\pi\,{}_3F_2\bigg(\begin{matrix} \tfrac12, \, \tfrac12, \, \tfrac12 \\ 1, \, \tfrac32 \end{matrix} \biggm|z\bigg)
=\int_0^1\!\int_0^1\frac{\d x\,\d y}{\sqrt{x(1-x)y(1-zxy)}},
\displaybreak[2]\\
\rho_1(z)
&=\pi\,{}_2F_1\bigg(\begin{matrix} \tfrac12, \, \tfrac12 \\ 1 \end{matrix} \biggm|z\bigg)
=\int_0^1\frac{\d x}{\sqrt{x(1-x)(1-zx)}},
\displaybreak[2]\\
\rho_2(z)
&=\pi\,{}_2F_1\bigg(\begin{matrix} -\tfrac12, \, \tfrac12 \\ 1 \end{matrix} \biggm|z\bigg)
=\int_0^1\frac{\sqrt{1-zx}}{\sqrt{x(1-x)}}\,\d x,
\end{align*}
then
$J_0(z)=\lambda(z)$,
\begin{align*}
J_1(z)&=-\frac{3+4z}{4z^2}\,\lambda-\frac{5(1-z)}{z^2}\,\rho_1+ \frac{13}{2z^2}\,\rho_2,
\\
J_2(z)&=\frac{105+480z+64z^2}{64z^4}\,\lambda+\frac{3151-2167z-984z^2}{144z^4}\,\rho_1-\frac{7247+3452z}{288z^4}\,\rho_2;
\end{align*}
in other words, each $J_n(z)$ is a $\mathbb Q(z)$-linear combination of $\lambda(z),\rho_1(z),\rho_2(z)$.
For $z^{-1}\in\mathbb Z\setminus\{\pm1\}$ we find out experimentally that the coefficients $a_n,b_n,c_n$ (depending, of course, on this~$z^{-1}$) in the representation
\[
J_n(z)=a_n\lambda(z)+b_n\rho_1(z)+c_n\rho_2(z)
\]
satisfy
\begin{equation*}
z^n2^{4n}a_n, \;  z^n2^{4n}D_{2n}^2b_n, \;  z^n2^{4n}D_{2n}^2c_n\in\mathbb Z
\qquad\text{for}\quad n=0,1,2,\dotsc.
\end{equation*}

Now observe that
\[
\det\begin{pmatrix} J_n & J_{n+1} \\ c_n & c_{n+1} \end{pmatrix}
=\det\begin{pmatrix} a_n & a_{n+1}  \\ c_n & c_{n+1} \end{pmatrix}\cdot\lambda(z)
+\det\begin{pmatrix} b_n & b_{n+1}  \\ c_n & c_{n+1} \end{pmatrix}\cdot\rho_1(z)
\]
for $n=0,1,2,\dots$\,. The sequences
\[
A_n=\det\begin{pmatrix} a_n & a_{n+1}  \\ c_n & c_{n+1} \end{pmatrix}
\quad\text{and}\quad
B_n=-\det\begin{pmatrix} b_n & b_{n+1} \\ c_n & c_{n+1} \end{pmatrix}
\]
satisfy the following third-order (again!) recurrence equation which is the exterior square of the recurrence for~$J_n$:
\begin{align*}
&
4(n+1)(n+2)^2(2n+1)^2(2n+3)^2z^8p_0(n)p_0(n-1) A_{n+1}
\displaybreak[2]\\ &\;
-4(n+1)^2(2n+1)^2z^4p_0(n-1)
\big(64(3z^2-20z+16)(27z-32)n^7
\\ &\;\quad
+ 64(3z^2-20z+16)(147z-170)n^6 + 16(3369z^3-26678z^2+44012z-20576)n^5
\\ &\;\quad
+ 16(2457z^3-20918z^2+34376z-15896)n^4
\\ &\;\quad
+ 4(843z^3-16808z^2+29432z-13736)n^3 - 4(1445z^3-6794z^2+9600z-4144)n^2
\\ &\;\quad
- (741z^3-6922z^2+10772z-4728)n + z^2(131z-66)\big)A_n
\displaybreak[2]\\ &\;
-n(2n-1)^2(1-z)z^2p_0(n+1)
\big(256(3z+8)(27z-32)n^8
\\ &\;\quad
- 256(3z+8)(15z-22)n^7 - 64(651z^2+661z-1744)n^6 + 192(59z^2-186)n^5
\\ &\;\quad
+ 16(1503z^2+697z-3610)n^4 - 16(79z^2-290z+116)n^3
\\ &\;\quad
- 4(569z^2-381z-580)n^2 + 4(11z^2-44z+18)n + 3(4z+3)(7z-10)\big)
A_{n-1}
\displaybreak[2]\\ &\;
-4(n-1)n^2(2n-3)^2(2n-1)^2(1-z)^2p_0(n)p_0(n+1) A_{n-2}
=0,
\end{align*}
where
\[
p_0(n)= 16(27z-32)n^4 + 48(13z-14)n^3 + 8(18z-11)n^2 - 4(19z-24)n - (7z+6)
\]
and
\begin{gather*}
A_0=\frac{13}{2z^2}, \quad
A_1=\frac{395z^2-1051z+591}{72z^6},
\\
A_2=\frac{15196z^4-201551z^3+548091z^2-543600z+183120}{3600z^{10}},
\\ \intertext{and}
B_0=0, \quad
B_1=\frac{1117z^2-2299z+1182}{72z^6},
\\
B_2=\frac{6867z^4-65547z^3+156430z^2-143530z+45780}{450z^{10}}.
\end{gather*}
Furthermore, by construction
\[
\lim_{n\to\infty}\frac{B_n}{A_n}=\frac{\lambda}{\rho_1}
\]
and, still only experimentally and for $z^{-1}\in\mathbb Z\setminus\{\pm1\}$,
\begin{equation*}
z^{2n+2}2^{2n}D_{2n}(n+1)(2n+1)^2A_n, \; z^{2n+2}2^{2n}D_{2n}^2(n+1)(2n+1)^2B_n\in\mathbb Z
\end{equation*}
for $n=0,1,2,\dots$\,.
In other words, the number $\lambda/\rho_1$ (but also the quotients $\lambda/\rho_2$ and $\rho_1/\rho_2$) are (weak) Ap\'ery limits for the values of $z$ in consideration.

\medskip
For real $k>0$ with $k^2\in\mathbb Z\setminus\{0,16\}$, the Mahler measure
\begin{align*}
\mu(k)
&=\mm(X+X^{-1}+Y+Y^{-1}+k)
\\
&=\frac1{(2\pi i)^2}\iint\limits_{|X|=|Y|=1}\log|X+X^{-1}+Y+Y^{-1}+k|\,\frac{\d X}{X}\,\frac{\d Y}{Y}
\end{align*}
is expected to be rationally proportional to the $L$-value
\[
L'(E,0)=\frac N{(2\pi)^2}\,L(E,2)
\]
of the elliptic curve $E=E_k:X+X^{-1}+Y+Y^{-1}+k=0$ of conductor $N=N_k=N(E_k)$.
This is actually proven \cite{BZ20} when $k=1$, $\sqrt2$, $2$, $2\sqrt2$ and $3$ for the corresponding elliptic curves \texttt{15a8}, \texttt{56a1}, \texttt{24a4}, \texttt{32a1} and \texttt{21a4} labeled in accordance with the database \cite{LMFDB};
the first number in the label indicates the conductor.

For the range $0<k<4$ we have the formula
\begin{equation*}
\mu(k)
=\frac k4\cdot{}_3F_2\biggl(\begin{matrix} \tfrac12, \, \tfrac12, \, \tfrac12 \\ 1, \, \tfrac32 \end{matrix}\biggm| \frac{k^2}{16} \biggr),
\end{equation*}
thus linking $\mu(k)$ to $z^{-1/2}\lambda(z)/\pi$ at $z=k^2/16$.
Furthermore, the quantity $z^{-1/2}\rho_1(z)$ in this case is rationally proportional to the imaginary part of the nonreal period of the same curve, while $z^{-1/2}\rho_2(z)$ is a $\mathbb Q$-linear combination of the imaginary parts of the nonreal period and the corresponding quasi-period.
It means that in many cases we can record $z^{-1/2}\rho_1(z)$ as a rational multiple of the central $L$-value of a quadratic twist of the curve $E$.
For example, when $k=2\sqrt2$ (hence $z=1/2$) the quadratic twist of the CM elliptic curve of conductor~32 coincides with itself and we have
\[
\lambda\Big(\frac12\Big)=2\sqrt2\,\pi L'(E,0)=16\sqrt2\,\frac{L(E,2)}{\pi}
\quad\text{and}\quad
\rho_1\Big(\frac12\Big)=4\sqrt2\,L(E,1),
\]
so that the recursion above with the choice $z=1/2$ realises the quotient $L(E,2)/\break(\pi L(E,1))$ as an Ap\'ery limit for an elliptic curve of conductor~32.
When $k=1$ we get
\[
\lambda\Big(\frac1{16}\Big)=8\pi L'(E,0)=30\,\frac{L(E,2)}{\pi}
\quad\text{and}\quad
\rho_1\Big(\frac1{16}\Big)=\frac12\,L(E,\chi_{-4},1)
\]
for the twist of the elliptic curve by the quadratic character $\chi_{-4}=\big(\frac{-4}{\cdot}\big)$;
this means that the quotient $L(E,2)/(\pi L(E,\chi_{-4},1))$ for an elliptic curve of conductor~15 is realised as an Ap\'ery limit.

Clearly, the range $0<k<4$ has a limited supply of elliptic $L$-values.
When $k>4$, one can write
\[
\mu(k)=\frac1{2\pi}\,f\bigg(\frac{16}{k^2}\bigg),
\]
where
\begin{align*}
f(z)
&=-\pi\bigg(\log\frac{z}{16}+\frac z4\,{}_4F_3\bigg(\begin{matrix} \tfrac32, \, \tfrac32, \, 1, \, 1  \\ 2, \, 2, \, 2 \end{matrix} \biggm|z\bigg) \bigg)
\\
&=-\int_0^1x^{-1/2}(1-x)^{-1/2}\log\frac{1-\sqrt{1-zx}}{1+\sqrt{1-zx}}\,\d x
\displaybreak[2]\\
&=\int_0^1\!\int_0^1\frac{x^{-1/2}(1-x)^{-1/2}(1-zx)^{1/2}y^{-1/2}}{1-(1-zx)y}\,\d x\,\d y
\\
&=Z\int_0^1\!\int_0^1\frac{x^{-1/2}(1-x)^{-1/2}(1-x/Z)^{1/2}(1-y)^{-1/2}}{x(1-y)+yZ}\,\d x\,\d y,
\end{align*}
with $Z=z^{-1}>1$.
At this point we see that the integrals resemble the integrals
\[
Z^{-l-m}\int_0^1\!\int_0^1\frac{x^j(1-x)^hy^k(1-y)^l}{(x(1-y)+yZ)^{j+k-m+1}}\,\d x\,\d y,
\]
with $h,j,k,l,m$ non-negative integers, appearing in the linear independence results for the dilogarithm \cite{RV05,VZ18}.
This similarity suggests looking at the family
\[
L_n(Z)=\int_0^1\!\int_0^1\frac{x^{n-1/2}(1-x)^{2n-1/2}(1-x/Z)^{1/2}y^n(1-y)^{n-1/2}}{(x(1-y)+yZ)^{n+1}}\,\d x\,\d y,
\]
where $Z=z^{-1}$ is a large (positive) integer. We tackle this double integral
by iterated applications of creative telescoping: while the first integration (no
matter whether one starts with $x$ or with $y$) can be done with the
Almkvist--Zeilberger algorithm, the second one requires more general
holonomic methods, since the integrand is not any more hyperexponential.
Using the \texttt{Mathematica} package \texttt{Holonomic\-Functions}~\cite{Kou10}, where these
algorithms are implemented, we find that the integral $L_n(Z)$ satisfies a
lengthy fourth-order recurrence equation. Moreover, it turns out that $L_n(Z)$
is a $\mathbb Q(Z)$-linear combination of $\rho_1=\rho_1(1/Z)$,
$\rho_2=\rho_2(1/Z)$, $\sigma_1=L_0(Z)$ and
\[
\sigma_2=\sigma_2(Z)
=\int_0^1\!\int_0^1\frac{x^{-1/2}(1-x)^{1/2}(1-x/Z)^{1/2}(1-y)^{1/2}}{x(1-y)+yZ}\,\d x\,\d y.
\]
One can produce a recurrence equation out of the one for $L_n(Z)$ to cast, for example, $\sigma_1/\rho_1$ as an Ap\'ery limit.
Because this finding does not meet any reasonable aesthetic requirements and does not imply anything (to be claimed) irrational, we leave it outside this note.

\end{document}